\documentclass[reqno, 12pt]{amsart}
\usepackage[body={7in,9.5in},top=1in, left=0.8in ]{geometry}
\usepackage{enumerate}
\usepackage{enumitem}
\usepackage{cancel}
\usepackage{enumitem}
\usepackage{amssymb}
\usepackage{mathdots}

\newtheorem{theorem}{Theorem}[section]
\newtheorem{definition}[theorem]{Definition}
\newtheorem{lemma}[theorem]{Lemma}
\newtheorem{proposition}[theorem]{Proposition}
\newtheorem{corollary}[theorem]{Corollary}
\newtheorem{example}[theorem]{Example}

\newcommand{\N}{\mathbb N}

\newcommand{\C}{\mathbb C}
\newcommand{\X}{\mathbb X}
\newcommand{\R}{\mathbb R}

\newcommand{\K}{\mathbb K}

\newcommand{\Z}{\mathbb Z}
\newcommand{\Lsc}{\mathcal{L}_s{(\mathbb C})}
\newcommand{\Lsci}{\mathcal{L}^1_s{(\mathbb C})}

\newcommand{\Lsk}{\mathcal{L}_s{(\mathbb K})}
\newcommand{\Lsku}{\mathcal{L}^u_s{(\mathbb K})}

\newcommand{\Fmultcu}{\X^1_0(\mathbb{C})}

\newcommand{\Fmultk}{\X_0(\mathbb{K})}

\newcommand{\Fc}{\mathbb{X}{(\mathbb C})}

\newcommand{\Fk}{\mathbb{X}{(\mathbb K})}

\DeclareMathOperator{\ord}{ord}

\newtheorem{df}[theorem]{Definition}

\theoremstyle{remark}
\newtheorem{remark}[theorem]{Remark}
\newtheorem{rmrk}[theorem]{Remark}

\numberwithin{equation}{section}

\begin{document}
\title[On Riordan groups involving formal semi-Laurent series]{On Riordan groups involving formal semi-Laurent series and their Lie group structure}

\subjclass[2020]{Primary: 13F25; Secondary: 13J05, 58B10, 22E66, 17B65.}

\keywords{Composition, composition inverse, formal power series, formal Laurent series, formal semi-Laurent series, Lie group, Lie algebra, Riordan array.}

\begin{abstract}

The main goal of this paper is to introduce and to investigate properties of generalized Riordan arrays and generalized Riordan groups that involve formal semi-Laurent series. In particular, we focus on the problem of isomorphy of generalized Riordan groups with the classical Riordan groups giving the negative answer to this problem. We also examine infinite dimensional Lie group structure of the group of generalized Riordan arrays.

\end{abstract}

\author[D.Bugajewski]{Dariusz Bugajewski}
\author[D.Bugajewski]{Dawid Bugajewski}
\author[X.-X.Gan]{Xiao-Xiong Gan}
\author[P.Ma\'ckowiak]{Piotr Ma\'ckowiak}

\address[D. Bugajewski, P. Ma\'ckowiak]{Department of Nonlinear Analysis and Applied Topology\\
  Faculty of Mathematics and Computer Science\\
  Adam Mickiewicz University\\
  Uniwersytetu Pozna\'nskiego 4\\
  61-614 Pozna\'n\\
  Poland}
	
\address[Dawid Bugajewski]{Faculty of Physics\\
University of Warsaw\\
Pasteura 5\\
02-093 Warsaw\\
Poland}
	
\address[X.-X.~Gan]{Department of Mathematics\\
Morgan State University\\
Baltimore, MD 21251\\
USA}

\email[D.~Bugajewski]{ddbb@amu.edu.pl}
\email[Dawid~Bugajewski]{d.bugajewski@student.uw.edu.pl}
\email[X.-X.~Gan]{gxiaoxio@gmail.com}
\email[P.~Ma\'ckowiak]{piotr.mackowiak@amu.edu.pl}

\maketitle

\section{Introduction}

The study of Riordan arrays, Riordan groups and their applications has been developed very extensively since the Riordan theory was established more than 40 years ago (see \cite{shwo}). More information concerning applications of Riordan arrays, especially in combinatorics, the interested reader can find, for example, in the monographs \cite{ba} and \cite{ssbchmw} or in the article \cite{cn} (see also the references given therein). It is also worth mentioning here the articles \cite{clmpm0,clmpm}, in which the Authors imposed a Fr\'echet - Lie group structure on the Riordan group and used it to obtain the matrix exponential of some lower triangular infinite matrices.

One can find several extensions of Riordan arrays and Riordan groups in the literature. Recall that the classical approach to Riordan arrays is based on the notion of formal power series and their composition. In this paper, we investigate generalized Riordan arrays and generalized Riordan groups that involve formal semi-Laurent series. Such types of generalizations of Riordan arrays from formal power series to certain wider class of formal semi-Laurent series or to formal Laurent series were investigated, for example, in the articles \cite{he11} and \cite{pmr}.

One of the main results of this article states that generalized Riordan groups are not isomorphic to the classical Riordan groups (see Section 6). We also show that some multiplicative groups of formal semi-Laurent series and formal power series are not isomorphic (see Section 4). These facts justify well the essence of generalizations of classical Riordan groups considered in this article. We also endow the set of generalized Riordan arrays with an infinite dimensional Lie group structure (over a locally convex space). We indicate its Lie algebra and its basic properties, including semisimplicity.

For completeness of the description of this paper let us add that in Section 3 we provide some comments concerning the composition of formal semi-Laurent series. In particular, we formulate some conditions for the existence of the composition of formal semi-Laurent series (following that proven in \cite{gank}) as well as some properties of such a composition like associativity or the Right Distributive Law. They will be useful for further investigation of Riordan arrays and Riordan groups.

\section{Preliminaries}

In this section we recall some basic definitions and facts which will be needed in the sequel. By $\Z$ we denote the set of integers, $\N$ is the set of positive integers, $\N_0:=\N\cup\{0\}$. The set of negative integers is denoted by $\Z_-$. If $n\in \N$, then $[n]:=\{1,\ldots,n\}$ and $[n]_0:=\{0,1,\ldots,n\}$. Let $\K\in\{\R,\C\}$. 

\begin{df}\label{def:101} (\cite{gb1}, \cite{g21}, \cite{hen})
A formal Laurent series $f$ over $\K$ is defined to be a mapping
from $\Z$ to $\K$. It is usually denoted by
\begin{equation}\label{eq:102}
\ldots + b_{-n} z^{-n} + \ldots + b_{-1} z^{-1} + b_0 + b_1 z + b_2 z^2 + \ldots
\end{equation}
or by $f(z) = \sum\limits_{n \in {\Z}} b_n z^n $, where $b_n \in {\K}$ for every $n \in \Z$. Also, we put $[z^n]f=f_n:=b_{n}$, $n\in\mathbb Z$. \\
We denote by ${\mathbb L}(\K)$, or simply by $\mathbb L$, the set of all formal Laurent series over $\K$. \\
A formal semi-Laurent series (briefly: fsLs) is a mapping $f:\Z\to \K$ such that $\#\{n\in \Z_{-}:\, f_n\neq 0\}<+\infty$. \\
Denote by ${\mathcal L}_s (\K)$ the set of all formal semi-Laurent series over $\K$.
\end{df}

By $0$ we briefly denote such a formal semi-Laurent series $g$ that $[z^n]g=0$ for all $n\in\mathbb Z$. Analogously, $1$ denotes such a $g\in\mathcal L_s$ that $[z^n]g=\delta_{n0}$. 
For each non-zero fsLs $g$ we associate with $g$ a special formal power series $t^g$ defined as $t^g_j=g_{\ord g+j},\,j\in \N_0$. Obviously, $\ord{t^g}=0$ and $g=z^{\ord g}t^g$.\\

It is clear that one can consider ${\mathbb X(\K)}$ (the set of all formal power series (briefly: fps) over $\K$) as a subset of ${\mathbb L}(\K)$ or ${\mathcal L}_s (\K)$. The basic algebraic and analytic structure of formal Laurent series was described in \cite{gb1}.
\\ It appears that $\mathcal L_s(\C)$ (and, by a simple analogy, $\mathcal L_s(\R)$) possesses a structure of a field.

\begin{theorem}\label{thm:202} (\cite{g21}, Th. 7.1.3, p. 192)
$\mathcal L_s(\C)$ is a field under addition and multiplication (the Cauchy product) and for a formal semi-Laurent series
$$
f(z) = f_m z^m + f_{m+1} z^{m+1} + f_{m+2} z^{m+2}  + \cdots,
$$
where $m \in \Z$ and $f_m \ne 0$, it holds
$$
f^{-1} (z) = z^{-m} (t^f)^{-1} (z)
$$
(recall that $t^f(z) = f_m + f_{m+1} z + f_{m+2} z^2 + \cdots$).
\end{theorem}

For any fsLs $f(z)$ its $n$-th power is denoted by $f^n(z)=\sum\limits_{j\in\mathbb Z}f^{(n)}_jz^j\in\mathcal L_s(\mathbb K)$, $n\in\N$. Based on the above theorem, we may admit the following: for every $n \in \N$ and $f \in \mathcal L_s (\K)$, we define $f^{-n} = (f^{-1})^n$.\\
Now, we extend the definition of the order of formal power series to formal Laurent series and formal semi-Laurent series.

\begin{df}\label{def:103} (cf. \cite{g21}, Def. 7.2.1, p. 199 and p. 233)
Let $f$ be a nonzero formal Laurent series. If $f \in {\mathcal L}_s (\K)$, then we define
$$
\ord f:=\min\{n\in \Z:\, f_n\neq 0\}.
$$
If $f \notin{\mathcal L}_s (\K)$, $\ord f := - \infty$. We also define $\ord 0:=+\infty$.
\end{df}

The following lemma describes a simple property of the order.

\begin{lemma}\label{pro:201}
Let $g, f \in \mathcal L_s (\K)$ be given such that
$$
g(z) = b_r z^r + b_{r+1} z^{r+1} + \cdots, \qquad f(z) = a_p z^p + a_{p+1} z^{p+1} + \cdots,
$$
where $b_r a_p \ne 0, \ r, p \in \Z$. Then
$$
\ord (gf) = r + p.
$$
\end{lemma}

\section{Composition of formal semi-Laurent series with formal power series}

The composition of formal Laurent series with formal power series was introduced in \cite{gb1} almost fifteen ago (some recent results concerning such a composition and its extension can be found in \cite{daw}). First, let us give the definition of composition of formal Laurent series with a formal power series, which extends that given in \cite{gb1} to the case of nonunit formal power series.

\begin{definition}\label{def:31}
Let $g(z) = \sum_{n\in\Bbb Z} b_n z^n$ be a formal Laurent series and let
\[
{\Bbb X}_g = \left\{ f(z) = \sum_{n=0}^\infty  a_n z^n  \in {\Bbb X}(\K) :  \sum_{n \in \Bbb Z}  b_n a_k^{(n)} \in \K \quad \mbox{for every }  k  \in \mathbb Z \right\} .
\]
If ${\Bbb X}_g \neq \emptyset $ (for example, if $g$ is a formal power series), we define $T_g : {\Bbb X}_g \longrightarrow \Bbb X (\mathbb K) $ by
\[
T_g(f)(z) = \sum_{k=-\infty}^\infty c_k z^k,
\]
where $c_k = \sum_{n \in \Bbb Z} b_n a_k^{(n)}  \in \K$ for every $k \in \mathbb Z$. $T_g(f)$ is called the composition of $g$ with $f$ and is also denoted by $g \circ f$.
\end{definition}

Notice that if $f$ is a unit (that is $a_0\neq 0$), then $f^n\in\mathbb X(\mathbb K)$ for every $n\in\mathbb Z$ and $g\circ f\in\mathbb X(\mathbb K)$ if it exists. The following result gives a necessary and sufficient condition for the existence of the composition of formal Laurent series with a unit formal power series.

\begin{theorem}\label{thm:201} (\cite{gb1})
Let $g(z) = \sum\limits_{k \in {\Z}}  b_k z^k$ be a formal Laurent series over~$\C$ and let $f(z) = \sum\limits_{k=0}^\infty a_k z^k$ be a formal power
series over $\C$ such that $ a_0 \not= 0$. Then the composition $g\circ f$ exists, if and only if both
$$
\sum_{n=k}^\infty \binom{n}{ k} b_{-n} a_0^{k-n} \in {\C} \quad \text{and} \quad  \sum_{n=k}^\infty \binom{n}{k} b_n a_0^{n-k} \in {\C}
$$
for every $k \in {\N} \cup \{0\}$.
\end{theorem}

\begin{corollary}\label{cor:201}   Let $g(z) = \sum\limits_{k=m}^\infty  b_k z^k$ be a formal
semi-Laurent series over $\C$ of order $m$ and let  $f(z) = \sum\limits_{k=0}^\infty a_k z^k$ be a formal power series. Then $(g\circ f)(z)\in {\mathcal L}_s (\C)$ if and only if
\begin{eqnarray}\label{cond1}
\sum_{n=k}^\infty \binom{n}{k} b_n a_0^{n-k} \in {\C}
\end{eqnarray}
for every $k \geq 0$.
\end{corollary}

\begin{proof}
If $a_0\neq 0$, then the claim follows immediately from Theorem \ref{thm:201}. Now let $a_0= 0$. Then (\ref{cond1}) holds, obviously. If $m\geq 0$, then $g\circ f$ exists since $g$ is a formal power series and $f$ is a nonunit. Suppose now that $m < 0$. Write $g = g^- + g^+ $, where
$$
g^-(z) = b_{-m} z^{-m} + b_{-m+1} z^{-m+1} + \ldots + b_{-1} z^{-1}, \quad g^+(z) = b_0 + b_1 z + b_2 z^2 + \ldots .
$$
Then $g^+\circ f\in\mathbb X(\mathbb C)$ and
$$
g^-\big(f(z)\big) = b_{-1} f^{-1}(z) + b_{-2} \big(f^{-1}(z)\big)^2 + \cdots + b_{-m} \big(f^{-1}(z)\big)^{m} \in \mathcal L_s.
$$
\end{proof}

\begin{rmrk}\label{cor:202}
Let $g(z) = \sum\limits_{k=m}^\infty  b_k \, z^k$ be a formal semi-Laurent series over $\C$. We showed that if $f$ is a nonunit, then $(g\circ f)(z) \in \mathcal L_s(\C)$. Also, notice that the coefficients of $g\circ f$ are then some polynomial functions of coefficients of $g$ and $f$.
\end{rmrk}

The following lemma describes simple properties of composition which will be useful in the sequel.

\begin{lemma}\label{lm:basic}
Let $g,h$ be fsLs and $f,k$ be nonunit formal power series. Then we have:
\begin{enumerate}[label={\textup{(\arabic*)}},ref=\textup{(\arabic*)}]
\item\label{bas:1} if $m,\, n\in \Z$ and $f\neq 0$, then $$(z^{m}\circ f)\cdot(z^{n}\circ f)=(z^{m+n}\circ f);$$
\item\label{bas:2} if $g\neq 0$ and $f\neq 0$, then
$$g\circ f=(z^{\ord g}\circ f)\cdot (t^g\circ f),$$
where $t^g_j:=g_{j+\ord g}$, $j\in \Z$;
\item\label{bas:3} (Right Distributive Law)\\
$$(g\circ f)\cdot(h\circ f)=(gh)\circ f;$$
\item\label{bas:4} if $g\neq 0$ and $f\neq 0$, and $m\in \Z$, then
$$(g\circ f)^m=g^m\circ f;$$
\item\label{bas:5}(Asssociativity of composition)\\
 $$(g\circ k)\circ f=g\circ (k\circ f).$$
\end{enumerate}
\end{lemma}

\begin{proof}
\ref{bas:1} It is clear that $z^m\circ f=f^m$ and $z^n\circ f=f^n$. Hence, $(z^m\circ f)\cdot (z^n\circ f)=f^mf^n=f^{m+n}=z^{m+n}\circ f$.

\ref{bas:2} Let $m=\ord g$. For any $j\in \Z$, we have
\begin{multline*}(g\circ f)_j=g_mf^{(m)}_j+g_{m+1}f^{(m+1)}_j+\ldots=g_mf^{(m)}_j+g_{m+1}f^{(m+1)}_j+\ldots+g_{j}f^{(j)}_j=\\(g_mf^{m}+g_{m+1}f^{m+1}+\ldots+g_{j}f^{j})_j=(f^{m}(g_m1+g_{m+1}f^{1}+\ldots+g_{j}f^{j-m}))_j,
\end{multline*}
so we get $g\circ f= f^{m}\cdot (t^g\circ f)=(z^m\circ f)\cdot (t^g\circ f)$.

\ref{bas:3} If $g=0$ or $h=0$, then the formula is true, obviously. Now, let $m=\ord g$ and $n=\ord h$. Since $g\circ f=(z^m\circ f)\cdot (t^g\circ f)$ and $h\circ f=(z^n\circ f)\cdot (t^h\circ f)$, we obtain \begin{multline*}(g\circ f)\cdot(h\circ f)=(z^m\circ f)\cdot (t^g\circ f)\cdot (z^n\circ f)\cdot (t^h\circ f)=
(z^m\circ f)\cdot (z^n\circ f)\cdot (t^g\circ f)\cdot (t^h\circ f)=\\(z^{m+n}\circ f)\cdot ((t^gt^h)\circ f)=(z^{m+n}(t^gt^h))\circ f=((z^{m}t^g)(z^{n}t^h))\circ f=(gh)\circ f,\end{multline*}
where we used the Right Distributive Law for formal power series (see \cite{bgm}, Section 4) (3rd equality) and property \ref{bas:2} (4th equality).

\ref{bas:4} By \ref{bas:3}, if $m>0$, $(g\circ f)^m=\underbrace{(g\circ f)\cdot \ldots \cdot (g\circ f)}_{m\times}=g^m\circ f$. If $m=0$, then the claim is true in a trivial way. Now, again by \ref{bas:3}, $(g\circ f)\cdot (g^{-1}\circ {f})=(gg^{-1})\circ f=1$, so $(g\circ f)^{-1}=g^{-1}\circ f$. Therefore, for $m<0$, $(g\circ f)^m=((g\circ f)^{-1})^{-m}=(g^{-1}\circ f)^{-m}=((g^{-1})^{-m})\circ f=g^{m}\circ f$.

\ref{bas:5} Let $m=\ord g$. If $m=+\infty$, the assertion is obviously true. Now, assume $m\in \Z$. We have\begin{multline*}
g\circ (k\circ f)=(z^{m}t^g)\circ (k\circ f)=(z^{m}\circ (k\circ f))\cdot(t^g\circ (k\circ f))=(k\circ f)^m\cdot ((t^g\circ k)\circ f)=\\
(k^m\circ f)\cdot((t^g\circ k)\circ f)=((z^m\circ k)\circ f)\cdot((t^g\circ k)\circ f)=((z^m\circ k)\cdot(t^g\circ k))\circ f)=\\((z^mt^g)\circ k)\circ f=(g\circ k)\circ f.
\end{multline*}
\end{proof}

\begin{lemma}\label{lem:201}
Let $g, f \in \mathcal L_s(\C)$ be given with $\ord (f) \geq 1$. Suppose that $g \circ f \in \mathcal L_s$. Then
\begin{equation}\label{eq:302}
\ord (g \circ f) = \ord g \cdot \ord f.
\end{equation}
\end{lemma}

\begin{proof} If $\ord g \geq 0$, then $g, f \in \X(\C)$ and the lemma is true (see \cite{g21}, Th. 1.2.5, p. 9). \\
Suppose that $\ord g = -m$ for some $m \in \N$. Since $\ord f \geq 1$, it follows that $\ord(f^{-1}) \leq -1$ (cf. Th. \ref{thm:202}) and then, denoting $g = g^- + g^+$ (see the proof of Corollary 3.3),
$$
\ord (g^-\circ f) = \ord (f^{-m}) = -m \ord f = \ord g \ord f.
$$
Thus (\ref{eq:302}) holds.
\end{proof}

The following example shows that the assumption $\ord (f) \geq 1$ in the above lemma is essential.
\begin{example} \label{ex:201}
Let $g(z) = \frac{1}{z} -2 + z, \ f(z) = 1 - z$. Then $\ord g = -1, \ord f = 0$. We have
$$
g\big(f(z)\big) = \frac{1}{1-z} -2 + (1 - z) = (1 + z + z^2 + \ldots ) -2 + 1 - z = z^2 + z^3 + \ldots .
$$
Then $\ord (g \circ f) = 2 \ne (-1) \cdot 0 = 0$.
\end{example}

\section{Multiplicative groups of formal semi-Laurent series and formal power series}

In this section we are going to show that some multiplicative groups of formal semi-Laurent series and formal power series are not isomorphic. These facts will show the essence of our generalizations which will be discussed in the next sections. \\

First, let $\Lsci:=\{f\in \Lsc:\, f_{\ord f}=1\}$ and $\Fmultcu:=\{f\in \Fc:\, f_0=1\}$.

\begin{theorem}\label{iso:1}
The groups $(\Lsci,\cdot)$ and $(\Fmultcu,\cdot)$, where $"\cdot"$ is the Cauchy product, are not isomorphic.
\end{theorem}

\begin{proof}
The group $\Lsci$ is isomorphic to $\Fmultcu\times\Z$ ($\mathbb Z$ -- the additive group of integers) via the isomorphism $\varphi:\Lsci\ni f\mapsto (z^{-\ord{f}}f,-\ord f)\in \Fmultcu\times \Z$, We also know that for any $f\in\Fmultcu$ and $m\in \N$, there exists $g\in \Fmultcu$ for which $f=g^m$ (see \cite[Cor. 9.2.9]{g21}), that is the group $\mathbb X_1(\mathbb C)$ is divisible. However, it is obvious that $\Fmultcu\times\Z$ is not divisible -- indeed, if $(f,n)\in\Fmultcu\times\Z$, $m\in\mathbb N$ and $m$ does not divide $n$, then there is no such $(g,k)\in \Fmultcu\times\Z$ that $(g,k)^m=(f,n)$, which completes the proof.
\end{proof}

Now, let us consider the following multiplicative groups $\Lsku:=\{f\in \Lsk:\, \ord{f}\in \Z\}$ and $\Fmultk:=\{f\in \Fk:\, \ord{f}=0\}$, where $\K\in \{\C,\R\}$.

\begin{theorem}\label{iso:2}
The multiplicative groups $\Lsku$ and $\Fmultk$ (endowed with Cauchy multiplication) are not isomorphic, $\K\in\{\R,\C\}$.
\end{theorem}
\begin{proof}
Let $f\in \Fmultk$. Then, for infinitely many $m\in \N$, there exist ${f}_0^{1/m}\in K$ which are $m$-th multiplicative roots of $f_0$ and, in consequence, for infinitely many $m\in \N$, there exist ${f}^{1/m}\in \Fmultk$ which are $m$-th multiplicative roots of $f$ (\cite[Cor. 9.2.9]{g21}). Now, suppose that $\varphi:\Fmultk\to \Lsku$ is an isomorphism and fix $F\in \Lsku$. There exists $f\in \Fmultk$ such that $F=\varphi(f)$ and $F=\varphi((f^{1/m})^m)=(\varphi(f^{1/m}))^m$ for infinitely many $m\in \N$. Hence, we get $\ord{F}=m\ord{\varphi(f^{1/m})}$, for infinitely many $m\in \N$, and the last equality is true for infinitely many $m\in \N$ only if $\ord{F}=0$. Therefore $\varphi(\Fmultk)\subset \Fmultk\neq \Lsku$, where we naturally embed $\Fmultk$ in $\Lsku$.
\end{proof}

\section{Generalized Riordan array space $RL$}

In this section we characterize Riordan arrays involving formal semi-Laurent series. First, let us provide a definition based on \cite{pmr}.
\begin{definition}\label{df:genR}
Let $g\neq 0$ be a fsLs and let $f$ be a fps such that $\ord f=1$. A Riordan Laurent matrix (RL matrix) or a generalized Riordan array associated with $g$ and $f$, denoted by $(g,f)$, is a mapping $(g,f): \mathbb Z\times  \mathbb Z\to \K$ defined as $$(g,f)_{m,n}:=[z^m](gf^n).$$
We call $g$ a generating fsLs or just a generating function of the generalized Riordan array $(g, f)$. The set of all RL matrices is
$$RL_1(\mathbb K):=\{(g,f):\,g, f\in\mathcal L_s,\,\,g\neq 0,\, \ord f=1\}.$$
Moreover, for $m\in \Z$, we denote
\[
RL_1^m(\mathbb K) := \{(g, f) \in RL_1 \ : \ \ord g = m\}.
\]
If $q=\ord{g}$, we can identify the RL matrix $(g,f)$ with a bi-infinite matrix whose columns are generated by fLs $ \ldots,gf^{-1},\,g,\,  gf,\, gf^2,\ldots$, that is
\begin{multline*}
(g,f)=\left[\begin{array}{llllllll}
\ddots&\vdots	&\vdots & \vdots&\vdots &\vdots &\vdots & \iddots \\
\ldots&	(gf^{-1})_{-1} &(gf^0)_{-1}& (gf^1)_{-1}& (gf^2)_{-1} & \ldots &(gf^m)_{-1}&\ldots  \\
\ldots&	(gf^{-1})_{0}&(gf^0)_{0}& (gf^1)_{0}& (gf^2)_{0} & \ldots &(gf^m)_{0}&\ldots  \\
\ldots&	(gf^{-1})_{1}&(gf^0)_{1}& (gf^1)_{1} & (gf^2)_{1} & \ldots &(gf^m)_{1}&\ldots\\
\ldots&	(gf^{-1})_{2}&(gf^0)_{2}& (gf^1)_{2} & (gf^2)_{2}  & \ldots &(gf^m)_{2}&\ldots\\	
\ldots&	\vdots&\vdots&  \vdots& \vdots & \ddots &\vdots& \vdots\\	
\ldots& (gf^{-1})_{m} &(gf^0)_{m}& (gf^1)_{m}& (gf^2)_{m}&\ldots & (gf^m)_{m}&\ldots \\
\iddots&\vdots	&\vdots&  \vdots& \vdots &  \vdots &\vdots& \vdots\\	\end{array}
\right]=\\\left[\begin{array}{llllllll}
\ddots&\vdots	&\vdots & \vdots&\vdots &\vdots &\vdots & \iddots \\
\ldots&(gf^{-1})_{q-1}	&0& 0& 0 & \ldots &0&\ldots  \\
\ldots&(gf^{-1})_{q+0}	&g_{q+0}& 0& 0 & \ldots &0&\ldots  \\
\ldots&(gf^{-1})_{q+1}	&g_{q+1}& (gf)_{q+1} & 0& \ldots &0&\ldots\\
\ldots&(gf^{-1})_{q+2}	&g_{q+2}& (gf)_{q+2} & (gf^2)_{q+2}  & \ldots &0&\ldots\\	
\ldots&\vdots	&\vdots&  \vdots& \vdots &  &\vdots& \vdots\\	
\ldots&\vdots	&\vdots&  \vdots& \vdots &  \ddots &0& \vdots\\	
\ldots&(gf^{-1})_{q+m}  &g_{q+m}& (gf)_{q+m}& (gf^2)_{q+m}&\ldots & (gf^m)_{q+m}&\ldots \\
\ldots&\vdots	&\vdots&  \vdots& \vdots &  \vdots &\vdots& \vdots\\	\end{array}
	\right],
\end{multline*}
because $(gf^n)_{q+m}=0$, for $n>m$, and $f^0:=1$.
\end{definition}

\begin{example}\label{ex:301} Let
$$
g(z) = z^{-3} + z^{-2} + z^{-1} + 1 + z + z^2 + \ldots, \quad f(z) = z + z^2 + \ldots\,.
$$
We have
$$
g(z) = z^{-3} (1 + z + z^2 +  \cdots ) = \frac{1}{z^3 (1-z)}, \quad f(z) = \frac{z}{1-z}.
$$
Then, applying the formal binomial series, we get
$$
g(z) f^j(z) = \frac{1}{z^3(1-z)} \cdot \left(\frac{z}{1-z}\right)^j = z^{j-3} (1-z)^{-(j+1)}=
$$

$$
= z^{j-3} \big[ 1 + (j+1) z + \frac{(j+1)(j+2)}{2!} z^2 + \ldots + \frac{(j+n)!}{j! n!} z^n + \ldots \big].
$$
for $j\in \mathbb Z$. Then
\begin{displaymath}
\left(\frac{1}{z^3 (1-z)}, \frac{z}{1-z}\right) = \left[ \begin{array}{ccccccccc}
  \ddots& \vdots &\vdots& \vdots& \vdots& \vdots& \vdots& \iddots\\
	  \ldots&1& 0& 0&0 &0 &0 &\ldots \\
  \ldots&-1& 1& 0&0 &0 &0 &\ldots \\
	\ldots&0& 0&1& 0& 0 & 0 & \ldots  \\
	\ldots&0& 0&1& 1 & 0 & 0 & \ldots  \\
	\ldots&0& 0&1& 2 & 1 & 0 & \ldots  \\	
	\ldots&0& 0&1& 3 & 3 & 1 & \ldots  \\	
   \iddots&\vdots &\vdots &\vdots & \vdots & \vdots & \vdots & \ddots  \end{array}
	\right].
\end{displaymath}
\end{example}

\begin{rmrk}\label{rmk:charact}
Sequential characterization of proper Riordan arrays is an important theoretical and practical issue. The characterization is as follows. An infinite lower triangular matrix $\hat{D}=[d_{i,j}]_{i,j\in \N_0}$ for which, $d_{i,i}\neq 0, i\in \N_0,$ is a proper Riordan matrix, that is, using our notation, $\hat{D}$ is the restriction of a generalized Riordan array $D=(g,f)\in RL_1^0(\mathbb K)$, for some fsLs $g$ and $f$ as in Definition \ref{df:genR}, to the set $\N_0\times \N_0$ if and only if there exist fps (sequences) $A=(a_0,a_1,a_2,\ldots)$, $a_0\neq 0$, and $Z=(z_0,z_1,z_2,\ldots)$ such that \begin{equation}\label{Acharact}d_{m+1,j+1}=a_0d_{m,j}+a_1d_{m,j+1}+a_2d_{m,j+2}+\ldots,\qquad m,j\in \N_0,\end{equation}
and
\begin{equation*}\label{Zcharact}d_{m+1,0}=z_0d_{m,0}+z_1d_{m,1}+z_2d_{m,2}+\ldots,\qquad m\in \N_0,\end{equation*}
(cf. Theorems 4.1, 4.4 and Remark 4.2 in \cite{ssbchmw}).

Let us put $D:=(g,f)\in RL_1(\mathbb K)$, for some given fsLs $g$ and $f$ as in Definition \ref{df:genR}, and fix $p\in \Z$.  Denote by $\hat{D}^p$ the restriction of $(t^{gf^p},f)\in RL_1^0(\mathbb K)$ to the set $\N_0\times \N_0$; $\hat{D}^p$ is a proper Riordan array associated with fsLs $t^{gf^p}$ and $f$. It is obvious that $gf^{p+n}=z^{\ord(gf^p)}t^{gf^p}f^n$, $n\in \N_0$. Hence, for $m, n\in \N_0$, we have
\begin{multline*}\label{relation}D_{m+\ord{(gf^p)},n+p}=(gf^{n+p})_{m+\ord{(gf^p)}}=(z^{\ord(gf^p)}t^{gf^p}f^n)_{m+\ord(gf^p)}=\\(z^{\ord(gf^p)})_{\ord(gf^p)}(t^{gf^p}f^n)_{m+\ord (gf^p)-\ord (gf^p)}=1\cdot(t^{gf^p}f^n)_{m}=\hat{D}^p_{m,n}.\end{multline*} 
Since $\hat{D}^p$ is a proper Riordan array, it is characterized by some sequences $A^p$ and $Z^p$ as described above (with $A^p$ and $Z^p$ in place of $A$ and $Z$, respectively).

Now, observe that if $k$, $k'$ are fsLs of order $0$ and $h$ is a fsLs of order $1$, then the Riordan arrays associated with $(k,h)$ and $(k',h)$ are characterized by exactly the same sequence $A$ for which equations like (\ref{Acharact}) hold, since $A$ depends solely on $h$ \cite{hespru} (it is independent of the fsLs $k$, the predecessor in the pair $(k,h)$ -- therefore the sequence $A^p$ is the same one (denote it as $A$) for all $p$). Further, if a column $q\in N_0$ of the array $(k,h)$ and a sequence $A$ that characterizes the array are known, then one can reconstruct the whole array $(k, h)$ with help of that $q$-th column and the sequence $A$ - this follows easily from the fact that $a_0\neq 0$ (see (\ref{Acharact})).

By the above remark and due to sequential characterization of the proper Riordan matrices (\cite[Theorems 2.1\&2.2]{hespru}) we arrive at the following theorem characterizing generalized Riordan arrays.

\begin{theorem}
Let $D=[d_{i,j}]_{i,j\in \Z}$ be a bi--infinite generalized lower diagonal matrix with non--zero diagonal whose entries belong to $\K$, that is, there exists $q\in \Z$ such that (diagonal elements) $d_{q+i,i}\neq 0,\, i\in \Z$, and $d_{q+i-j,i}=0,\, i\in \Z,\,j\in \N$. The matrix $D$ is a generalized Riordan array, that is, $D=(g,f)\in RL_1(\mathbb K)$ for some fsLs $g$, $\ord g=q$, and $f$, $\ord{f}=1$, if and only if there exist sequences $A=(a_0,a_1,\ldots),$ $a_0\neq 0$, and $Z^0=(z_0,z_1,\ldots)$ in $\K$ such that 
\begin{equation}\label{Acharactgen}d_{m+1,j+1}=a_0d_{m,j}+a_1d_{m,j+1}+a_2d_{m,j+2}+\ldots,\qquad m,j\in \Z,\end{equation}
and
\begin{equation}\label{Zcharactgen}d_{q+m+1,0}=z_0d_{q+m,0}+z_1d_{q+m,1}+z_2d_{q+m,2}+\ldots,\qquad m\in \Z\setminus\{-1\}.\end{equation}
\end{theorem}
\end{rmrk}
Let us notice that in the last formula there is characterized $0$-th column of $D$. Without loss of generality we could take any column of $D$, say $q$-th, and characterize it with some sequence $Z^q$ (possibly different from $Z^0$). 

Finally, observe that to determine whether a bi--infinite generalized lower triangular array $D$ (with non--zero elements on its diagonal) is a generalized Riordan array it is necessary and sufficient to know a column of it and a sequence $A$ for which the condition (\ref{Acharactgen}) is satisfied or, equivalently, to know the diagonal element of a column, say $q$-th column for some $q\in \Z$, and sequences $A$ together with $Z^q$ for which the conditions (\ref{Acharactgen}) and (\ref{Zcharactgen}) are satisfied (with a slight adaptation of indices in (\ref{Zcharactgen})).

In the sequel, however, we will use the following equivalent characterization of the Riordan arrays, independent of the matrix representation.
\begin{definition}\label{df:product}
A Riordan Laurent matrix (or: generalized Riordan array) is defined as a pair $(g,f)$, where $g\in\mathcal L_s$, $g\neq 0$, $f\in\mathbb X(\mathbb C)$, $\ord f=1$. \\
Let $(g,f)$ and $(h,k)$ be two RL matrices associated with $g,f$ and $h,k$, respectively. The product of $(g,f)$ and $(h,k)$ is the RL matrix defined as follows
$$
(g,f) * (h,k) := (g\cdot(h\circ f),k\circ f).
$$
\end{definition}
Sometimes we will omit the symbol $'*'$ and write the product of $(g,f)$ and $(h,k)$ simply as $(g,f)(h,k)$. \\

\begin{rmrk}\label{rmk4.1}
\begin{enumerate}
\item Note that the compositions $h\circ f$ and $k\circ f$ are well-defined, since $\ord f=1$, and we have $\ord (g\cdot(h\circ f))=\ord g+\ord (h\circ f)=\ord g+\ord h$ and $\ord k\circ f =1$ (see Lemma \ref{lem:201} and Proposition \ref{pro:201}). Thus, the above definition of the product $*:RL_1(\mathbb K)\times RL_1(\mathbb K)\to RL_1(\mathbb K)$ is correct. Moreover, if $(g, f) \in RL^{m_1}_1(\mathbb K)$ and $(h, l) \in RL^{m_2}_1(\mathbb K)$, then
$$
(g, f) \ast (h, l) \in RL^{m_1 + m_2}_1(\mathbb K).
$$
\item Theorem 19(b) from \cite{pmr} implies that the above definition is equivalent to Definition \ref{df:genR}.
\end{enumerate}
\end{rmrk}
\section{Generalized Riordan groups}

Let us mention that $RL_1^0(\mathbb K)$ may be identified with the well-known classical Riordan group (see e.g. \cite{shwo} and the comment at the beginning of Remark \ref{rmk:charact}). In a similar way, one can prove the following

\begin{theorem}\label{thm:RLgroup}
The set $RL_1(\K)$ endowed with the binary operation of multiplication $*$, defined in Definition \ref{df:product}, is a group. The identity element of the group is the RL matrix $(1,z)$ and, for $(g,f)\in RL_1(\K)$, the inverse of $(g,f)$, denoted by $(g,f)^{-1}$, is defined by the equality
$$(g,f)^{-1}=(g^{-1}\circ (f^{[-1]}), f^{[-1]}),$$
where the formal power series  $f^{[-1]}$ is the compositional inverse of $f$, that is, $f\circ f^{[-1]}=f^{[-1]}\circ f=z$.
\end{theorem}

Now, we are going to focus on the problem of isomorphicity of some groups considered in this paper. Let us denote by $\mathbb X_1^z(\mathbb C)$ -- the group of all nonunit formal power series satisfying $[z^1]g\neq 0$ with composition $\circ:(g_1,g_2)\mapsto g_2\circ g_1$ (this seemingly reverse order will be useful here). We also denote, for $f\in\mathbb X(\mathbb C)$ with $\ord(f)=1$, $\frac{z}{f}=(z^{-1}f)^{-1}$.

\begin{theorem}\label{lemR1}
$RL_1(\mathbb C)$ is isomorphic to the semidirect product $(\mathbb X_0(\mathbb C)\times \mathbb Z)\rtimes_{\Psi}\mathbb X_1^z(\mathbb C)$, where the homomorphism $\Psi:\mathbb X_1^z(\mathbb C)\rightarrow \mathrm{Aut}(\mathbb X_0(\mathbb C)\times \mathbb Z)$ is given by
\[
\Psi(f)(g,n)=((g\circ f)\cdot \left(\frac{z}{f}\right)^n,n).
\]
\end{theorem}

\begin{proof}
Let us first prove that for every $f\in\mathbb X_1^z(\mathbb C)$, $\Psi(f)\in\textrm{Aut}(\mathbb X_0(\mathbb C)\times \mathbb Z)$ and that $\Psi$ is a homomorphism. It is easy to check that every mapping $\Psi(f)$ is a bijection; moreover, by the Right Distributive Law,
\[
\Psi(f)((g_1,n_1)(g_2,n_2))=\Psi(f)(g_1g_2,n_1+n_2)=(((g_1g_2)\circ f)\cdot \left(\frac{z}{f}\right)^{n_1+n_2},n_1+n_2)=
\]
\[
((g_1\circ f)\cdot \left(\frac{z}{f}\right)^{n_1},n_1)((g_2\circ f)\cdot \left(\frac{z}{f}\right)^{n_2},n_2)=\Psi(f)(g_1,n_1)\Psi(f)(g_2,n_2).
\]
Also, for every $f_1,f_2\in\mathbb X_1^z(\C)$, $g\in\mathbb X_0(\C)$, $n\in\mathbb Z$,
\[
\Psi(f_2)\Psi(f_1)(g,n)=((((g\circ f_1)\circ f_2)\cdot \left(\frac{f_2}{f_1\circ f_2}\right)^n)\cdot \left(\frac{z}{f_2}\right)^n,n)=\Psi(f_1\circ f_2)(g,n).
\]
Notice that this actually means that $\Psi$ is a homomorphism (not an antihomomorphism), because the group action in $\mathbb X_1^z(\mathbb C)$ was defined as $(g_1,g_2)\mapsto g_2\circ g_1$.\\
Now, the postulated isomorphism can be given by $\chi:RL_1\ni (g,f)\mapsto ((t^g,-\textrm{ord}(g)),f)\in(\mathbb X_0(\mathbb C)\times \mathbb Z)\rtimes_{\Psi}\mathbb X_1^z(\mathbb C)$ (recall that $t^g=z^{-\mathrm{ord}(g)}g$): indeed, $\chi$ is obviously bijective and ($n:=-\mathrm{ord}(g)$, $m:=-\mathrm{ord}(h)$)
\[
\chi((g,f))\chi((h,k))=((t^g,n),f)((t^h,m),k)=
\]
\[
((t^g\cdot(t^h\circ f)\cdot \left(\frac{z}{f}\right)^m,n+m),k\circ f)=
\chi((g\cdot (h\circ f),k\circ f))=\chi((g,f)(h,k)).
\]
\end{proof}

Let us now denote by $\mathbb X_1^1(\mathbb C)$ -- the group of all nonunit formal power series with $f_1=1$ and composition as the group action. Further, let
$$RL_{1,1}(\mathbb K):=\{(g,f)\in RL_1: f_1=1\}.$$
and
\[
RL_{1,1}^0(\mathbb K) := \{(g, f) \in RL_1^0 \ :  f_1=1\}.
\]

\begin{corollary}\label{cor:62}
Arguing similarly as in the proof of Theorem \ref{lemR1} one can prove that $RL_{1,1}(\mathbb C)$ is isomorphic to the semidirect product $(\mathbb X_0(\mathbb C)\times \mathbb Z)\rtimes_{\Psi}\mathbb X_1^1(\mathbb C)$.
\end{corollary}

\begin{lemma}\label{lm:power}
Let $(g,f)\in RL_1(\mathbb K)$ and $m\in \N$. Then $(g,f)^m=(\prod_{i=0}^{m-1} g\circ f^{[i]},f^{[m]})$, where $f^{[0]}:=z$ and $f^{[i]}:=\underbrace{f\circ\ldots\circ f}_{i \text{ copies of }f}$, $i\geq 1$.
\end{lemma}

\begin{proof}
It is clear that the formula is valid for $m=1$. Suppose that it works for $m=1,\ldots, k-1$. Then, for $m=k$, we get
\begin{multline*}(g,f)^k=(g,f)(g,f)^{k-1}=(g,f)(\prod_{i=0}^{k-2} g\circ f^{[i]},f^{[k-1]})=(g((\prod_{i=0}^{k-2} g\circ f^{[i]})\circ f), f^{[k-1]}\circ f)=\\(g(\prod_{i=0}^{k-2} g\circ f^{[i+1]}), f^{[k-1]}\circ f)=(\prod_{i=0}^{k-1} g\circ f^{[i]}, f^{[k]}),
\end{multline*}
where the penultimate equality is due to the Right Distributive Law.
\end{proof}

\begin{lemma}\label{lemR2}
The group $RL^0_{1,1}(\mathbb C)$ is divisible.
\end{lemma}

\begin{proof}
Let $(g,f)\in RL^0_{1,1}(\mathbb C)$ and $n\in\mathbb Z_+$. The solution $(G,F)^n=(g,f)$ with unknowns $G,F$ can be written as
\[
(G\cdot (G\cdot (G\cdot\ldots\circ F)\circ F)\circ F,F^{[n]})=(g,f).
\]
A formal power series $F$ satisfying $F^{[n]}=g$ exists (see e.g. \cite{g21}, Thm. 8.3.5, Prop. 8.1.10); we will now prove that for such a $F$, an appropriate $G$ exists. By Lemma \ref{lm:power}, the analyzed equation can be written  in the form
\[
G\cdot (G\circ F)\cdot\ldots\cdot(G\circ F^{[n-1]})=g;
\]
comparing the coefficients of both sides, we obtain: $G_0^n=g_0$ ($[z^0]$), which obviously possesses a solution in $\mathbb C$, $G_1=\frac{g_1}{nG_0^{n-1}}$ ($[z^1]$) and, for every $k>1$,
\[
nG_0^{n-1}G_k+\Omega_k(G_0,...,G_{k-1})=g_k
\]
($\Omega_k:\mathbb C^k\rightarrow \mathbb C$ -- some polynomials). Therefore we obtain $G$ by a recursive calculation of its coefficients.
\end{proof}

\begin{corollary}\label{cor:noniso}
$RL_{1,1}(\mathbb C)$ and $RL^0_{1,1}(\mathbb C)$ are not isomorphic.
\end{corollary}

\begin{proof}
$RL_{1,1}(\mathbb C)$ is not divisible by Corollary \ref{cor:62} and the fact that $\mathbb Z$ is not. The claim is now a direct result of Lemma \ref{lemR2}.
\end{proof}

\begin{lemma}\label{lm:decomp}
Let $(g,f)\in RL_1(\mathbb K)$ and $m\in \N$. Then $(g^m,f)=(g,z)^m(1,f)$.
\end{lemma}
\begin{proof}
It is obvious that $(g,z)^m=(g^m,z)$. Hence, $(g,z)^m(1,f)=(g^m,z)(1,f)=(g^m(1\circ z),f\circ z)=(g^m,f)$.
\end{proof}

\begin{lemma}\label{lm:homprop1}
Let $\varphi=(\varphi_1,\varphi_2):RL_1^0(\mathbb K)\to RL_1(\mathbb K)$ be a homomorphism. Then, for every $k\in \K\setminus\{0\}$, $\ord\varphi_1(k,z)=\ord\varphi_1(1,kz)=0$.
\end{lemma}

\begin{proof}
Let us fix $k\in \K\setminus\{0\}$. It is clear that for infinitely many $m\in\N$ there exist $k^{1/m}\in \K$, which are $m$-th roots of $k$. In what follows we consider only such $m$'s.
Now, set $a_m:=\varphi_1(k^{1/m},z),\, b_m:=\varphi_2(k^{1/m},z)$. Since $(k^{1/m},z)^m=(k,z)$, by Lemma \ref{lm:power} we obtain $(\varphi_1(k,z),\varphi_2(k,z))=\varphi(k,z)=(\varphi(k^{1/m},z))^m=(a_m,b_m)^m=(\prod_{i=0}^m a_m\circ b_m^{[i]},b_m^{[m]})$.
The fact that $\ord{b_m}=1$ and Lemma \ref{lem:201} imply that $\ord{\varphi_1(k,z)}=\ord{\prod_{i=0}^m a_m\circ b_m^{[i]}}=\ord{a_m}+\ord{a_m\circ b_m}+\ord{a_m\circ b_m^{[2]}}+\ldots+\ord{a_m\circ b_m^{[m
-1]}}=m\ord a_m$ for all the $m$'s under consideration. However, this is possible only if $\ord \varphi_1(k,z)=0$.

Now, let us investigate $\ord \varphi_1(1,kz)$. Observe that $$\varphi(1,kz)=\varphi(1,\underbrace{k^{1/m}z\circ \ldots \circ k^{1/m}z}_{m\times})=\varphi((1,k^{1/m}z)^m)=(\varphi(1,k^{1/m}z))^m$$ and, using similar reasoning as in the previous part of the proof, we get $\ord\varphi_1(1,kz)=m\ord \varphi_1(1,k^{1/m}z)$ for all relevant $m$'s, which results in $\ord\varphi_1(1,kz)=0$.
\end{proof}

\begin{lemma}\label{lm:homprop2}
Let $\varphi=(\varphi_1,\varphi_2):RL_1^0(\mathbb K)\to RL_1(\mathbb K)$ be a homomorphism. Then $\ord\varphi_1(1,f)=\ord\varphi_1(1,\frac{1}{f_1}f)$, $(1,f)\in RL_1^0(\mathbb K)$.
\end{lemma}

\begin{proof}
We have
$$
\varphi\left(1,\frac{1}{f_1}f\right)= \varphi\left(1,\left(\frac{1}{f_1}z\right)\circ f\right)=\varphi\left((1,f)\left(1, \frac{1}{f_1}z\right)\right)=\varphi(1,f)\varphi\left(1, \frac{1}{f_1}z\right).
$$
Hence, $\ord \varphi_1\left(1,\frac{1}{f_1}f\right)=\ord \varphi_1(1,f)+\ord \varphi_1\left(1,\frac{1}{f_1}z\right)$ and, by Lemma \ref{lm:homprop1}, the assertion is true.
\end{proof}

\begin{lemma}\label{lm:itercomp} Let $f\in \X(\K)$, $\ord{f}=1$ and $f_1=1$. Then there exists an iterative square root $g\in\X(\K)$ of $f$, that is, $f=g\circ g$ with $g_1=1$ and $\ord{g}=1$. \end{lemma}
\begin{proof} We shall construct a fps $g$ by induction. Let $g_0:=0,\, g_1:=1$. Observe that whatever the values of $g_n,\, n\geq 2$ are, we have $f_0=(g\circ g)_0$ and $f_1=(g\circ g)_1$. Now, suppose that we have determined $g_0,\ldots, g_n$ so that to have $f_i=(g\circ g)_i$, $i\in [n]_0$. It holds that $f_{n+1}=(g\circ g)_{n+1}$ if and only if $f_{n+1}=\sum_{i=1}^{n+1}g_i (g^i)_{n+1}$. The last equality is equivalent to $f_{n+1}=g_1g_{n+1}+g_{n+1}g_1^{n+1}+\sum_{i=2}^{n}g_i (g^i)_{n+1}=2g_{n+1}+\sum_{i=2}^{n}g_i (g^i)_{n+1}$, where $\sum_{i=2}^{n}g_i (g^i)_{n+1}$ is independent of $g_{n+1}$. From this it follows that for $g_{n+1}:=(f_{n+1}-\sum_{i=2}^{n}g_i (g^i)_{n+1})/2$ we have $f_{n+1}=(g\circ g)_{n+1}$, which completes the proof. \end{proof}

\begin{lemma}\label{lm:homprop3}
Let $\varphi=(\varphi_1,\varphi_2):RL_1^0(\K)\to RL_1(\K)$ be a homomorphism. Then $\ord\varphi_1(g,f)=\ord\varphi_1(1,f)$ for all $(g,f)\in RL_1^0(\K)$.
\end{lemma}

\begin{proof}
For any $(g,f)\in RL_1^0(\K)$ denote $\varphi(g,f)=(\varphi_1(g,f),\varphi_2(g,f))$, where $\varphi_1(g,f)\in\Lsk\setminus\{0\},\, \varphi_2(g,f)\in\Lsk$; it is clear that $\ord{\varphi_2(g,f)}=1$.

Let us fix $(g,f)\in RL_1^0(\K)$. Since $g=g_0\hat{g}$, where $\hat{g}_i:=g_i/g_0,\,i\in \N$, by Lemma \ref{lm:decomp}, we get
\begin{multline*} \varphi(g,f)=\varphi((g,z)(1,f))=\varphi(g,z)\varphi(1,f)=\varphi(g_0\hat{g},z)\varphi(1,f)=\\\varphi(g_01,z)\varphi(\hat{g},z)\varphi(1,f)=\varphi(g_01,z)(\varphi(\hat{g}^{1/m},z))^m\varphi(1,f)=(a_0,b_0)(a_m,b_m)^m(c,d)=(\star),
\end{multline*}
where $(a_0,b_0):=\varphi(g_0,z), (a_m,b_m):=\varphi(\hat{g}^{1/m},z)$, $m\in \N$, $(c,d):=\varphi(1,f)$, and $\hat{g}^{1/m}$ is an $m$-th multiplicative root of $\hat{g}$ (which is known to exist, since $\hat{g}_0=1$). By Lemma \ref{lm:power}, we obtain
\begin{multline*}\left(\star\right)=\left(a_0,b_0\right)\left(\prod_{i=0}^{m-1} a_m\circ b_m^{[i]},b_m^{[m]}\right)\left(c,d\right)=\left(a_0\left(\left(\prod_{i=0}^{m-1} a_m\circ b_m^{[i]}\right)\circ b_0\right),b_m^{[m]}\circ b_0\right)\left(c,d\right)=\\\left(a_0\left(\left(\prod_{i=0}^{m-1} a_m\circ b_m^{[i]}\right)\circ b_0\right)\left(c\circ\left(b_m^{[m]}\circ b_0\right)\right),d\circ\left(b_m^{[m]}\circ b_0\right)\right)=\\\left(a_0\left(\prod_{i=0}^{m-1} a_m\circ\left(b_m^{[i]}\circ b_0\right)\right)\left(c\circ\left(b_m^{[m]}\circ b_0\right)\right),d\circ\left(b_m^{[m]}\circ b_0\right)\right),
\end{multline*}
where we used the Right Distributive Law to obtain the last equality. Now, by Proposition \ref{pro:201}, Lemma \ref{lem:201} and Lemma \ref{lm:homprop1}, for all $m\in \N$, we get
\begin{multline*}
\ord \varphi_1\left(g,f\right)=\ord{a_0\left(\prod_{i=0}^{m-1}a_m\circ\left(b_m^{[i]}\circ b_0\right)\right)\left(c\circ\left(b_m^{[m]}\circ b_0\right)\right)}=\\\ord a_0+\ord{\left(\prod_{i=0}^{m-1} a_m\circ\left(b_m^{[i]}\circ b_0\right)\right)}+\ord{\left(c\circ\left(b_m^{[m]}\circ b_0\right)\right)}=\ord a_0+m\ord a_m+\ord c=\\\ord\varphi_1\left(g_01,z\right)+m\ord\varphi_1\left(\hat{g}^{1/m},z\right)+\ord \varphi_1\left(1,f\right)=m\ord\varphi_1\left(\hat{g}^{1/m},z\right)+\ord \varphi_1\left(1,f\right).
\end{multline*}
Therefore we obtain
\begin{equation}\label{eqn:ord}
\ord \varphi_1(g,f)=\ord \varphi_1(1,f), \text{ for all }(g,f)\in RL_1^0(\K).
\end{equation}
\end{proof}

\begin{theorem}\label{iso:3}
Suppose that $\varphi=(\varphi_1,\varphi_2):RL_1^0(\K)\to RL_1(\K)$ is a homomorphism. Then $\varphi(RL_1^0(\K))\subset RL_1^0(\K)$.
\end{theorem}

\begin{proof}
Suppose that there exists $(g,f)\in RL_1^0(\K)$ for which $\ord \varphi_1(g,f)\neq 0$. By Lemma \ref{lm:homprop3}, we have $0\neq\ord\varphi_1(g,f)=\ord\varphi_1(1,f)$ and, by Lemma \ref{lm:homprop2}, we can assume that $f_1=1$. By Lemma \ref{lm:itercomp}, it follows that $f$ possesses $2^m$-th iterative roots $f^{[1/2^m]}$ with $f^{[1/2^m]}_0=1$, $m\in \N$. For any $m\in \N$, we have $$\varphi(1,f)=\varphi(1,\underbrace{f^{[1/2^m]}\circ\ldots\circ f^{[1/2^m]}}_{2^m\times})=\varphi((1,f^{[1/2^m]})^{2^m})=(\varphi(1,f^{[1/2^m]}))^{2^m}.$$
Therefore, by Lemmas \ref{lm:power} and \ref{lm:homprop1}, Proposition \ref{pro:201} and Lemma \ref{lem:201} we obtain (cf. the proof of Lemma \ref{lm:homprop3})
$$\ord \varphi_1(1,f)=2^m\ord \varphi_1(1,f^{[1/2^m]}),\quad m\in \N,$$
which implies that $\ord \varphi_1(1,f)=0$, contradicting our assumption. The proof is complete.

\end{proof}

From the above theorem we obtain the following

\begin{corollary}
The multiplicative group $RL_1(\K)$ is not isomorphic to the Riordan group $RL_1^0(\K)$.
\end{corollary}

From the proofs of Lemma \ref{lm:homprop3} and Theorem \ref{iso:3} we conclude that the next corollary is true.

\begin{corollary}\label{cor:noniso1}
Suppose that $A\subset RL_1^0(\K)$ is a subgroup such that whenever $(g,f)\in A$, then $(g,z)\in A,\,(1,f)\in A$, for infinitely many $m\in \N$: $(g^{1/m},z)\in A$, and for infinitely  many $n\in \N$: $(1,f^{[1/n]})\in A$. If $\varphi:A\to RL_1(\K)$ is a homomorphism, then $\varphi(A)\subset RL_1^0(\K)$.
\end{corollary}

Notice that Corollary \ref{cor:noniso1} implies Corollary \ref{cor:noniso}.

\section{Infinite dimensional Lie group structure of the generalized Riordan group}

In the following paragraph, we are going to endow the set of generalized Riordan arrays with an infinite dimensional Lie group structure (over a Fr\'{e}chet or a locally convex space). For the sake of clarity, we are often going to simply write {\it Lie group} instead of {\it infinite dimensional/Fr\'{e}chet/locally convex Lie group} if it does not cause any ambiguity.

Let $s$ denote the Fr\'{e}chet space of all complex valued sequences and let $s_L$ denote the space of all mappings $a:\mathbb Z\rightarrow \mathbb C$ such that there exists a $N_{\ord a}\in\mathbb Z$ satisfying $a(n)=0$ if $n<N_{\ord a}$. Obviously $s_L$ is a locally convex topological vector space, with a separating family of seminorms $||a||_m=|a(m)|$, $m\in\mathbb Z$. Let us emphasize that $s_L$ is not a Fr\'{e}chet space -- if it was, the sequence $(a_k)_{k\in\mathbb N}$, where $a_k(m)=1$ if $m\geq -k$ and $0$ if $m<-k$ would be a Cauchy sequence, but it has no limit in $s_L$ (it would require us to consider the broader space $\mathbb L$ of {\it all} formal Laurent series, which is however, not a group with formal multiplication since $gf$ is not well-defined for all $g,f\in\mathbb L$ (see e.g. \cite{gb1})). Fortunately it is enough to model an infinite dimensional manifold (see e.g. \cite{monastir}) and, consequently, a Lie group. \\

Let us emphasize that the constructed objects (e.g. the Lie algebra of the generalized Riordan group) will possess a somewhat predictable (due to the results from previous sections) semidirect product structure. Nevertheless, for the convenience of the reader and to establish a firm base for further investigations, we will present all steps of explicit calculations (including e.g. full proofs of continuity and differentiability of the composition of fsLs with fps), dropping some possible shortcuts resulting from the semidirect product structure.

Let us introduce the following

\begin{proposition}\label{lie1}
$((RL_1(\mathbb C),*,\tau),A)$, where $*$ -- the multiplication of Riordan arrays, $\tau$ -- topology on $RL_1(\mathbb C)$ equivalent to the product topology on $s_L\times s$, $A$ -- maximal atlas compatible with $\left\{(RL_1(\mathbb C),\varphi)\right\}$, $\varphi:(g,f)\mapsto((...,0,g_{\mathrm{ord}(g)},g_{\mathrm{ord}(g)+1},\ldots,g_0,g_1,\ldots),(f_1,f_2,f_3,...))$, is an infinite dimensional Lie group modeled on the locally convex topological vector space $s_L\times s$.
\end{proposition}
\begin{proof}
Obviously $\varphi(RL_1(\mathbb C))$ is an open subset of $s_L\times s$ and $\varphi$ is a homeomorphism onto its image. It is also easy to check the smoothness of mappings $*$ and $(g,f)\mapsto (g,f)^{-1}$: for every $(g,f),(h,k)\in RL_1(\mathbb C)$, $(g,f)(h,k)=(g\cdot(h\circ f),k\circ f)$ and $(g,f)^{-1}=\left(\frac{1}{g\big(f^{[-1]}\big)}, f^{[-1]} \right)$; see that all operations on formal power series/semi-Laurent series in these equations are smooth (holomorphic -- they involve only polynomial/rational functions of the series' coefficients, see e.g. \cite{daw2}, Lemma 3.3).
\end{proof}
Let us also notice that a sequence of Riordan arrays $((g_n,f_n))_{n\in\mathbb N}$, $(g_n,f_n)\equiv (a_{ij}^{(n)})_{i,j}$ is convergent with respect to topology $\tau$, if and only for all $i,j$, the sequence $(a_{ij}^{(n)})_{n\in\mathbb N}$ is convergent. This is because for all $i,j$ and $n\in\mathbb N$, $a_{ij}^{(n)}$ is a polynomial function of the coefficients of $g_n$ and $f_n$.\\
In what follows, we are not going to distinguish formal series $f$ and their sequence representations $\varphi(f)$ (if it does not cause any ambiguity) to simplify the notation.
\begin{lemma}\label{diffofcomp}
Let $g\in\mathcal L_s$ and $h,\chi$ be nonunit formal power series. Then
\begin{enumerate}
\item $\lim\limits_{t\rightarrow 0}g\circ(h+t\chi)=g\circ h$,
\item $\lim\limits_{t\rightarrow 0}\frac{g\circ(h+t\chi)-g\circ h}{t}=(g'\circ h)\chi$.
\end{enumerate}
\end{lemma}
\begin{proof}
We will first prove (2); (1) then is its obvious consequence. It follows immediately from the definition of fsLs multiplication and the formula for multiplicative inverse of fps (cf. \cite{g21}, Thm. 1.1.8 with an obvious observation $\left(z^n(a_0+a_1z+...)\right)^{-1}=z^{-n}(a_0+a_1z+...)^{-1}$) that multiplication and taking inverse in $\mathcal L_s$ are continuous operations with respect to the topology under consideration. Moreover, all sums over $s$ below are finite due to the nonunitness assumptions. Therefore
\[
\lim\limits_{t\rightarrow 0}\frac{g\circ (h+t\chi)-g\circ h}{t}=\sum\limits_{s\in\mathbb Z}g_s\left(\lim\limits_{t\rightarrow 0}\frac{(h+t\chi)^s- h^s}{t}\right).
\]
If $s=0$, then obviously $(h+t\chi)^s- h^s=1-1=0$ for all $t$; now, if $s\geq 1$, then
\[
\frac{(h+t\chi)^s- h^s}{t}=\sum\limits_{k=0}^{s-1}\binom{s}{k}t^{s-k-1}h^k\chi^{s-k}\rightarrow sh^{s-1}\chi\mbox{ as }t\rightarrow 0.
\]
Also, as $t\rightarrow 0$,
\[
\frac{(h+t\chi)^{-1}-h^{-1}}{t}=\frac{hh^{-1}(h+t\chi)^{-1}-(h+t\chi)(h+t\chi)^{-1}h^{-1}}{t}=-\chi h^{-1}(h+t\chi)^{-1}\rightarrow -\chi h^{-2}
\]
and for $s<-1$, by the mathematical induction, we can deduce that
\[
\frac{(h+t\chi)^s-h^s}{t}=h^{-1}\frac{(h+t\chi)^{s+1}-h^{s+1}}{t}-\frac{\chi}{h}(h+t\chi)^s\rightarrow s\chi h^{s-1}.
\]
As a result,
\[
\lim\limits_{t\rightarrow 0}\frac{g\circ (h+t\chi)-g\circ h}{t}=\sum\limits_{s\in\mathbb Z}g_ss\chi h^{s-1}=\left(\sum\limits_{s\in\mathbb Z}g'_sh^s\right)\chi=(g'\circ h)\chi,
\]
which completes the proof.
\end{proof}
\begin{lemma}
Let $*$ denote the multiplication in $RL_1(\mathbb C)$. Then
\begin{enumerate}
\item for every $(g,f)\in RL_1(\mathbb C)$,
\[
T_{(g,f)}RL_1(\mathbb C)=\left\{(\alpha,\beta):\alpha\in\mathcal L_s(\mathbb C),\beta\in\mathbb X(\mathbb C),[z^0]\beta=0\right\}:=\mathcal L_s(\mathbb C)\times z\mathbb X(\mathbb C).
\]
\item the tangent mapping $T*:TRL_1(\mathbb C)\times TRL_1(\mathbb C)\rightarrow TRL_1(\mathbb C)$ is given by
\begin{eqnarray*}\label{tan}
&&T*(((f,h),(\phi,\chi)),((g,k),(\gamma,\kappa)))=\\&&((f\cdot(g\circ h),k\circ h),(\phi\cdot (g\circ h)+f\cdot (\gamma\circ h)+f\cdot(g'\circ h)\cdot\chi,(k'\circ h)\cdot \chi+\kappa\circ h))
\end{eqnarray*}
for every $(f,h),(g,k)\in RL_1(\mathbb C)$, $(\phi,\chi)\in T_{(f,h)}RL_1(\mathbb C)$, $(\gamma,\kappa)\in T_{(g,k)}RL_1(\mathbb C)$.
\end{enumerate}
\end{lemma}
\begin{proof}
The first part of the claim is obvious. Now, using Lemma \ref{diffofcomp}, we have
\begin{eqnarray*}
&&
\lim\limits_{t\rightarrow 0}\frac{(f+t\phi,h+t\chi)(g+t\gamma,k+t\kappa)-(f,h)(g,k)}{t}=\\
&&
\lim\limits_{t\rightarrow 0}\frac{((f+t\phi)\cdot(g+t\gamma)\circ (h+t\chi),(k+t\kappa)\circ (h+t\chi))-(f\cdot(g\circ h),k\circ h)}{t}=\\
&&
\lim\limits_{t\rightarrow 0}\left(\frac{f\cdot(g+t\gamma)\circ(h+t\chi)-f\cdot(g\circ h)}{t}+\phi\cdot(g+t\gamma)\circ(h+t\chi)\right.,\\
&&
\left.\frac{k\circ(h+t\chi)-k\circ h}{t}+\kappa\circ(h+t\chi)\right)=\\
&&
\left(\phi\cdot(g\circ h)+f\cdot(\gamma\circ h)+f\cdot(g'\circ h)\chi,(k'\circ h)\chi+\kappa\circ h \right)\in T_{(f,h)(g,k)}RL_1.
\end{eqnarray*}
\end{proof}
Now we use the notations from \cite{monastir} to construct the Lie algebra corresponding to the defined Lie group $RL_1(\mathbb C)$. For every $(f,h)\in RL_1(\mathbb C)$, $((g,k),(\gamma,\kappa))\in TRL_1(\mathbb C)$,
\[
(f,h).((g,k),(\gamma,\kappa))=T*(((f,h),(0,0)),((g,k),(\gamma,\kappa)))=((f\cdot(g\circ h),k\circ h),(f\cdot(\gamma\circ h),\kappa\circ h)).
\]
Therefore the left-invariant vector field corresponding to some $(\gamma,\kappa)\in T_{(1,z)}RL_1(\mathbb C)$ can be expressed as $X_{(\gamma,\kappa)}:RL_1\ni(f,h)\mapsto ((f,h),\tilde{X}_{(\gamma,\kappa)}(f,h))\in TRL_1(\mathbb C)$, where
\[
X_{(\gamma,\kappa)}(f,h)=(f,h).((1,z),(\gamma,\kappa))=((f,h),\underbrace{(f\cdot(\gamma\circ h),\kappa\circ h)}_{\tilde{X}_{(\gamma,\kappa)}(f,h)}).
\]
Denoting $X_{1,2}:=X_{(\gamma_{1,2},\kappa_{1,2})}$, we have for instance
\begin{eqnarray*}
&&
(d_{\tilde{X}_1}\tilde{X}_2)(1,z)=\lim\limits_{t\rightarrow 0}\frac{\tilde{X}_2((1,z)+(t\tilde{X}_1(1,z)))-\tilde{X}_2((1,z))}{t}=\\
&&
\lim\limits_{t\rightarrow 0}\frac{((1+t\gamma_1)\cdot(\gamma_2\circ(z+t\kappa_1)),\kappa_2\circ(z+t\kappa_1))-(\gamma_2,\kappa_2)}{t}=\\
&&
\lim\limits_{t\rightarrow 0}\left(\frac{\gamma_2\circ(z+t\kappa_1)-\gamma_2}{t}+\gamma_1\cdot(\gamma_2\circ(z+t\kappa_1)),\frac{\kappa_2\circ(z+t\kappa_1)-\kappa_2}{t}\right)=(\gamma_2'\kappa_1+\gamma_1\gamma_2,\kappa_2'\kappa_1).
\end{eqnarray*}
This leads us to the following
\begin{proposition}
The Lie algebra $rl_1(\mathbb C)$ corresponding to the Lie group $RL_1(\mathbb C)$ is $\mathcal L_s(\mathbb C)\times z\mathbb X(\mathbb C)$ with Lie bracket
\[
[(\gamma_1,\kappa_1),(\gamma_2,\kappa_2)]=(\gamma_2'\kappa_1-\gamma_1'\kappa_2,\kappa_2'\kappa_1-\kappa_1'\kappa_2).
\]
\end{proposition}
The following proposition establishes some basic properties of $rl_1(\mathbb C)$.
\begin{proposition}\label{liealg}
\begin{enumerate}
\item[(i)] $rl_1(\mathbb C)$ is an infinite-dimensional semisimple (but not simple) Lie algebra;
\item[(ii)] $rl_1(\mathbb C)=l_s(\mathbb C)\rtimes_\psi x_0^z(\mathbb C)$, where $l_s(\mathbb C)$ is the Lie algebra of all formal semi-Laurent series with trivial Lie bracket, $x_1^z(\mathbb C)$ is the Lie algebra corresponding to the Lie group $\mathbb X_1^z(\mathbb C)$, that is the vector space of all nonunit formal power series with Lie bracket
    \[
    [f_1,f_2]=f_2'f_1-f_2f_1',
    \]
    and $\psi:x_1^z(\mathbb C)\rightarrow \textrm{Der}(l_s(\mathbb C))$ is a homomorphism defined by
    \[
    \psi_{\kappa}(\gamma)=\kappa\gamma'
    \] 
($\textrm{Der}(l_s(\mathbb C))$ denotes the set of derivations on the Lie algebra	$l_s(\mathbb C)$).  	
\end{enumerate}
\end{proposition}
\begin{proof}
To prove $rl_1(\mathbb C)$ is not simple, it is enough to observe that e.g. $l_s(\mathbb C)\times \left\{f\in\mathbb X(\mathbb C):\ord f\geq N\right\}$ for every $N\geq 2$ is its proper ideal due to an inequality $\ord(\kappa_2'\kappa_1-\kappa_1'\kappa_2)\geq \max(\ord\kappa_1,\ord\kappa_2)$. Let us now assume $I$ to be an abelian proper ideal of $rl_1(\mathbb C)$. Then for every $(\gamma_1,\kappa_{1}),(\gamma_2,\kappa_{2})\in I$, $\kappa_2'\kappa_1-\kappa_1'\kappa_2=0$. Let now $\kappa=\kappa_2\kappa_1^{-1}\in\mathcal L_s$. See that $\kappa_1\kappa_1^{-1}=1$ implies $(\kappa_1^{-1})'=-\kappa_1'\kappa_1^{-2}$, and therefore $\kappa'=(\kappa_2'\kappa_1-\kappa_1'\kappa_2)\kappa_1^{-2}=0$, which means $\kappa_2=c\kappa_1$ for some $c\in\mathbb C$. Now, $[rl_1(\mathbb C),I]\subseteq I$ implies that for every $\kappa_2\in x_1^z(\mathbb C)$ and $\kappa_1\in I$, $\kappa_1,\kappa_2\neq 0$, $\ord \kappa_2>\ord\kappa_1$, there exists such $c\in\mathbb C$ that $\kappa_2'\kappa_1-\kappa_1'\kappa_2=c\kappa_1$. Then $\ord (c\kappa_1)=\ord (\kappa_2'\kappa_1-\kappa_1'\kappa_2)>\ord \kappa_1$, which implies $c=0$, $\kappa_2=\tilde{c}\kappa_1$ for some $\tilde{c}\in\mathbb C\setminus\{0\}$, which contradicts the assumption $\ord \kappa_2>\ord\kappa_1$. This proves (i). Claim (ii) is a direct consequence of the definition of semi-direct product of Lie algebras.
\end{proof}
\begin{remark}
The Lie algebra corresponding to the classical Riordan group (with fps instead of fsLs) was expressed e.g. in \cite{clmpm0}, Prop. 22 as the space of linear mappings of the form $L_{\chi,\alpha}:\mathbb X(\mathbb C)\ni h \mapsto \chi h+\alpha h'$, $\chi,\alpha\in\mathbb X(\mathbb C)$, $\alpha_0=0$. See that 
\begin{eqnarray*}
&&
[L_{\alpha_1,\chi_1},L_{\alpha_2,\chi_2}]h=\chi_1(\chi_2h+\alpha_2h')+\alpha_1(\chi_2h+\alpha_2h')'-\chi_2(\chi_1h+\alpha_1h')-\alpha_2(\chi_1h+\alpha_1h')'\\
&&
=(\chi_2'\alpha_1-\chi_1'\alpha_2)h+(\alpha_1\alpha_2'-\alpha_2\alpha_1')h'=L_{\chi_2'\alpha_1-\chi_1'\alpha_2,\alpha_1\alpha_2'-\alpha_2\alpha_1'}h.
\end{eqnarray*}
This shows that the Lie algebra of generalized Riordan matrices proposed in this paper is an extension of the earlier constructed Lie algebra corresponding to the classical Fr\'{e}chet-Lie Riordan group.
\end{remark}
\begin{remark}
It is easy to endow the group $Rl_{1,1}(\mathbb C)$ with a Lie group structure as well - one only has to change the map $\varphi$ to $\varphi:(g,f)\mapsto((...,0,g_{\mathrm{ord}(g)},g_{\mathrm{ord}(g)+1},\ldots,g_0,g_1,\ldots),(f_2,f_3,...))$. The corresponding Lie algebra is 
\[
rl_{1,1}(\mathbb C)=\left\{(\alpha,\beta):\alpha\in\mathcal L_s(\mathbb C),\beta\in\mathbb X(\mathbb C),[z^0]\beta=[z^1]\beta=0\right\}:=\mathcal L_s(\mathbb C)\times z^2\mathbb X(\mathbb C).
\]
with the Lie bracket analogous to $rl_1(\mathbb C)$ and properties analogous to those from Proposition \ref{liealg}.
\end{remark}


\begin{thebibliography}{10}

\bibitem{ba}
Barry, P. {\it Riordan arrays: a prime}, Kildare, Ireland: Logic Press, 2016.

\bibitem{daw}
Bugajewski, Dawid, {\it The inverse and the composition in the set of formal Laurent series}, arXiv:2202.13948v2.

\bibitem{daw2}
Bugajewski, Dawid, {\it Topology and geometry of the general composition of formal power series - towards Fr\'echet-Lie group - like formalism}, arXiv:2409.09853v3.

\bibitem{bgm}
Bugajewski, D., Galimberti, A. and P. Ma\'ckowiak, On composition and Right Distributive Law for formal power series of multiple variables, Aequat. Math., {\bf 99}, 2025, 21--35.

\bibitem{cn}
Cameron, N. and Nkwanta, A., {\it Riordan matrices and lattice path enumeration},  Not. Am. Math. Soc, {\bf 70}(2), 2023, 231--243.

\bibitem{clmpm0} Cheon, G. S., Luzón, A., Morón, M. A., Prieto-Martinez, L. F. and Song, M., {\it Finite and infinite dimensional Lie group structures on Riordan groups}, Adv. in Math., {\bf 319}, 2017, 522--566.

\bibitem{clmpm}
Chocano, P.J., Luz\'on A., Mor\'on, M.A. and Prieto - Martinez, L.F., {\it Characteristic curves and the exponentation in the Riordan Lie group: A connection through examples}, J. Math. Anal. Appl. {\bf 532}, 2024, 127989, 1 - 18.


\bibitem{g21}
Gan, X., {\it Formal Analysis, An introduction}, De Gruyter, 2021, ISBN: 978-3-110-59785-1.

\bibitem{gb1}
Gan, X. and Bugajewski, D., {\it On formal Laurent series}, Bull. Braz. Math. Soc., {\bf 42}(3), 2011, 415--437.

\bibitem{gank}
Gan, X. and Knox, N., {\it On composition of formal power series}, Int. J. Math. and Math. Sci. {\bf 30}(12), 2002, 761-770.

\bibitem{he11}
He, T.-X., {\it Riordan arrays associated with Laurent series and generalized Sheffer-type groups},
Linear Alg. and its Appl., {bf 435}, 2011, 1241-1256. 

\bibitem{hespru}
He, T.-X., Sprugnoli, R., {\it Sequence characterization of Riordan arrays}, Discrete Math., {\bf 309}(12), 2009, 3962--3974.

\bibitem{hen}
Henrici, P., {\it Applied and Computational Complex Analysis}, John Wiley and Sons, 1988.

\bibitem{monastir}
Neeb, K.–H., {\it Monastir Summer School: Infinite-Dimensional Lie Groups}, Monastir, 2005.

\bibitem{pmr}
Pieto - Martinez, L.F. and Ricob, J., {\it Bi-infinite Riordan matrices: a matricial approach to multiplication and composition of Laurent series}
arXiv:2504.07593v1.

\bibitem{shwo}
Shapiro, L., Getu, S., Woan, W. and Woodson, L., {\it The Riordan group},
Discrete. Appl. Math., {\bf 34}, 1991, 229--239.
		
\bibitem{ssbchmw}		
Shapiro, L., Sprugnoli, R., Barry, P., Cheon, GS., He, TX., Merlini, D., and W. Wang, {\it The Riordan Group and Applications}, Springer, 2022.

\end{thebibliography}
\end{document}